\documentclass[11pt]{article}
\usepackage{latexsym}
\usepackage{amsmath}
\usepackage{amsthm,color}
\usepackage{amssymb}
\usepackage{mathrsfs}
\usepackage[colorlinks,  linkcolor=blue,  anchorcolor=blue, citecolor=blue]{hyperref}

\usepackage{lscape}
\newtheorem{theorem}{\indent Theorem}[section]

\newtheorem{definition}[theorem]{\indent Definition}
\newtheorem{lemma}[theorem]{\indent Lemma}

\newtheorem{case*}[theorem]{\indent Case*}
\begin{document}
\renewcommand{\baselinestretch}{1.3}


\begin{center}
    {\large \bf The Nehari manifold for fractional p-Laplacian system involving concave-convex nonlinearities and sign-changing weight functions}
\vspace{0.5cm}\\
{\sc Maoding Zhen  }\\
{\small 1) School of Mathematics and Statistics, Huazhong University of Science and Technology,\\ Wuhan 430074, China}\\
{\small 2) Hubei Key Laboratory of Engineering Modeling and Scientific Computing, Huazhong University of Science and Technology,}
{\small Wuhan, 430074, China}\\
\end{center}


\renewcommand{\theequation}{\arabic{section}.\arabic{equation}}
\numberwithin{equation}{section}


\begin{abstract}
In this paper, we consider a fractional p-Laplacian system \eqref{int1} with both concave-convex nonlinearities and sign-changing weight functions in bounded domains. With the help of the Nehari\ manifold, we prove that the system has at least two nontrivial solutions when the pair of the parameters $(\lambda,\mu)$ belongs to a certain subset of $\mathbb{R}^{n}$ .

\end{abstract}
\vspace{-1 cm}

\footnote[0]{ \hspace*{-7.4mm}
$^{*}$ Corresponding author.
The authors were supported by the NSFC grant 11571125\\
E-mails:d201677010@hust.edu.cn }

\section{Introduction}
In this paper, we study the following system involving fractional p-Laplacian:
\begin{equation}\label{int1}
\begin{cases}
(-\Delta)_{p}^{s}u= \lambda f(x)|u|^{q-2}u+\frac{2\alpha}{\alpha+\beta}h(x)|u|^{\alpha-2}u|v|^{\beta}  & \text{in} \ \Omega,\\
(-\Delta)_{p}^{s}v= \mu g(x)|v|^{q-2}v+\frac{2\beta}{\alpha+\beta}h(x)|u|^{\alpha}|v|^{\beta-2}v&\text{in} \ \Omega,\\

u=v=0&\text{in}\ \mathbb{R}^{n}\setminus\Omega,
\end{cases}
\end{equation}
where $\Omega$ is a smooth bounded domain in $\mathbb{R}^{n}$, $s\in(0,1),\ n>ps$, $1<q<p$ and $\alpha>1, \beta>1$ satisfy $p<\alpha+\beta<p^{\star}=\frac{np}{n-ps}$, where $p^{\star}$ is the fractional Sobolev exponent and $(-\Delta)_{p}^{s}$ is the fractional p-Laplacian operator which is defined as $$(-\Delta)_{p}^{s}u(x)=2\lim\limits_{\epsilon\rightarrow0}\int_{\mathbb{R}^{n}\setminus B_{\epsilon}(x)}\frac{|u(y)-u(x)|^{p-2}(u(y)-u(x))}{|x-y|^{n+ps}}dy,\ \ \ x\in\mathbb{R}^{n}.$$  The pair of parameters $(\lambda,\mu)\in\mathbb{R}^{n}\setminus(0,0)$ and the weight functions $f,g,h$ satisfy the following conditions;\\
(A) $f,g\in L^{q^{\star}}(\Omega)$ where $q^{\star}=\frac{\alpha+\beta}{\alpha+\beta-q}$ and $f^{+}=\max\{\pm f,0\}\neq 0$ or $g^{+}=\max\{\pm g,0\}\neq 0$;
(B) $h\in C(\overline{\Omega})$ with $\|h\|_{\infty}=1$ and $h\geq 0$.

Recently, a great deal of attention has been focused on studying of equations or systems involving fractional Laplacian and corresponding nonlocal problems, both for their interesting theoretical structure and their concrete applications(see \cite{BCPS,CDS2,MF,CS1,CRS,CS2,ZW} and references therein). This type of operator arises in a quite natural way in many different contexts, such as, the thin obstacle problem, finance, phase transitions, anomalous diffusion, flame propagation and many others(see\cite{DPV,MMC,LS} and references therein).

Compared to the Laplacian  problem, the fractional p-Laplacian problem is nonlocal and more challenging. On the one hand, for the fractional elliptic problems when $p=2$ has been investigated by many researchs. For example, C. Br\"{a}ndle, E. Colorado, A. de Pablo and U. S\'{a}nchez \cite{CCPS} studied the fractional Laplacian equation involving concave-convex nonlinearity for the subcritical case, they prove that there exists a finite parameter $\Lambda>0$ such that for $0<\lambda<\Lambda$ there exist at least two solutions, for $\lambda=\Lambda $ there exists at least one solution and for $\lambda>\Lambda$ there is no solution. Furthermore, B. Barrios, E. Colorado, A. de Pablo and U. S\'{a}nchez \cite{BCPS} studied the nonhomogeneous equation involving fractional Laplacian and proved the existence and multiplicity of solutions under suitable conditions of $s$ and $q$. E. Colorado, A. de Pablo and U. S\'{a}nchez \cite{CDS2} studied the fractional equation with critical Sobolev exponent, they proved that the existence and the multiplicity of solutions under appropriate conditions on the size of $f$. For more other advances on this topic, see \cite{SV11,SV12} for the subcritical, \cite{SV13,SV14} for the critical case.

For the fractional p-Laplacian equation, S. Goyal, K. Sreenadh \cite{SK} studied the following equation involving concave-convex nonlinearities and sign-changing weight functions.
\begin{equation*}
\begin{cases}
(-\Delta)_{p}^{s}u= \lambda h(x)  |u|^{q-1}u+b(x)|u|^{r-1}u $$  &\text{in} \ \ \Omega,\\

u=0& \text{on} \ \mathbb{R}^{n}\setminus\Omega,\\
\end{cases}
\end{equation*}
they showed that the existence and multiplicity of solutions by
minimization on the suitable subset of Nehari manifold using the fibering maps and proved that there exists $\lambda_{0}$ such that for $\lambda\in(0,\lambda_{0})$, it has at least two non-negative solutions.

B. Cheng, X. Tang \cite{CXT} studied the existence of solutions for the following fractional p-Laplacian equation with sign-changing potential and nonlinearity
\begin{equation*}
(-\Delta)_{p}^{s}u+v(x)|u|^{p-2}u=f(x,u)\ \ \forall x\in\mathbb{R}^{N},
\end{equation*}
 where $p\geq2,\  N\geq2,\ 0<s<1,\ V\in C(\mathbb{R}^{N},\mathbb{R})$ and $f\in C(\mathbb{R}^{N}\times\mathbb{R},\mathbb{R})$, under the suitable conditions they prove the equation has at least one nontrivial solution.

It is also natural to study the coupled system of equations. W. Chen, S. Deng \cite{WSD} considered the special case when $f(x)=g(x)=h(x)=1$ for system \eqref{int1}. They prove that the system admits at least two nontrivial solutions under proper conditions of $\lambda$ and $\mu$.

The purpose of this paper is to study system \eqref{int1} in the case of $2<\alpha+\beta<p^{\star}$, by variational methods and a Nehari manifold decomposition, we prove that the system admits at least two nontrivial solutions when the pair of parameters $(\lambda,\mu)$ belongs to certain subset of $\mathbb{R}^{2}$.

To express the main results, we introduce $$\Theta=\{z\in\mathbb{R}^{2}\setminus{(0,0)}\mid 0<(|\lambda|\|f\|_{L^{q^{\star}}})^{\frac{p}{p-q}}+(|\mu|\|g\|_{L^{q^{\star}}})^{\frac{p}{p-q}}<C(\alpha,\beta,p,q,S)\}$$
and $C(\alpha,\beta,p,q,S)=\left[\frac{p-q}{2(\alpha+\beta-q)}S^{\frac{\alpha+\beta}{p}}\right]^{\frac{p}{\alpha+\beta-p}}\left(S^{-\frac{q}{p}}\frac{\alpha+\beta-q}{\alpha+\beta-p}\right)^{-\frac{p}{p-q}}$.\\
$$\Psi=\{z\in\mathbb{R}^{2}\setminus{(0,0)}\mid 0<(|\lambda|\|f\|_{L^{q^{\star}}})^{\frac{p}{p-q}}+(|\mu|\|g\|_{L^{q^{\star}}})^{\frac{p}{p-q}}<D(\alpha,\beta,p,q,S)\}$$
and $D(\alpha,\beta,p,q,S)=\left(\frac{q}{p}\right)^{\frac{p}{p-q}}\left[\frac{p-q}{2(\alpha+\beta-q)}S^{\frac{\alpha+\beta}{p}}\right]^{\frac{p}{\alpha+\beta-p}}\left(S^{-\frac{q}{p}}\frac{\alpha+\beta-q}{\alpha+\beta-p}\right)^{-\frac{p}{p-q}}.$\\

Where $S$ is the best Sobolev constant that will be introduced later.

Our main results are:
\begin{theorem}\label{TH1}
Suppose that the weight functions $f,g,h$ be satisfied with the conditions $(A)$ and $(B)$, for each $(\lambda,\mu)\in \Theta$, then system \eqref{int1} has at least one nontrivial solution in $W^{s,p}(\Omega)\times W^{s,p}(\Omega)$.
\end{theorem}
\begin{theorem}\label{TH2}
Suppose that the weight functions $f,g,h$ be satisfied with the conditions $(A)$ and $(B)$, for each $(\lambda,\mu)\in \Psi$, then system \eqref{int1} has at least two nontrivial solutions in $W^{s,p}(\Omega)\times W^{s,p}(\Omega)$.
\end{theorem}

The paper is organized as follows. In section \ref{sec2}, we introduce some preliminaries. In section \ref{sec3}, we define the Nehari manifold and give some Lemmas that will be used later. In section \ref{sec4}, we prove the existence of Palais-Smale sequence. In section \ref{sec5}, we give the results of local minimization problem for system \eqref{int1}. Finally, the proofs of Theorem \ref{TH1} and Theorem \ref{TH2} are given in section \ref{sec6}.

\section{ Preliminaries }\label{sec2}
In this section, we introduce some preliminaries that will be used to establish the energy functional for system \eqref{int1}.
Let $s\in(0,1)$ and $p\in[1,+\infty)$, we define the usual fractional Sobolev space $ W^{s,p}(\Omega)$ endowed with the norm
\begin{align}\label{P1}
\|u\|_{W^{s,p}(\Omega)}=\|u\|_{L^{p}(\Omega)}+\left(\int_{\Omega\times\Omega}\frac{|u(x)-u(y)|^{p}}{|x-y|^{n+ps}}dxdy\right)^{\frac{1}{p}}.
\end{align}
Let $Q=\mathbb{R}^{2n}\setminus(\mathcal{C}_{\Omega}\times\mathcal{C}_{\Omega})$ with $\mathcal{C}_{\Omega}=\mathbb{R}^{n}\setminus\Omega$. We define $$X=\left\{u\mid u:\mathbb{R}^{n}\rightarrow \mathbb{R} \ is \ measurable \ and \ \int_{Q}\frac{|u(x)-u(y)|^{p}}{|x-y|^{n+ps}}dxdy<+\infty \right\}.$$
The space $X$ is endowed with the norm defined by
\begin{align}\label{P2}
\|u\|_{X}=\|u\|_{L^{p}(\Omega)}+\left(\int_{Q}\frac{|u(x)-u(y)|^{p}}{|x-y|^{n+ps}}dxdy\right)^{\frac{1}{p}}.
\end{align}
The functional space $X_{0}$ denotes the closure of $C^{\infty}_{0}(\Omega)$ in $X$. By Theorem 6.5 and Theorem 7.1 in \cite{DPV} the space $X_{0}$ is a Hilbert space which can be endowed with the scalar product defined for any $\phi,\psi\in X_{0}$ as
\begin{align}\label{P4}
\langle \phi,\psi\rangle_{X_{0}}=\int_{Q}\frac{|\phi(x)-\phi(y)|^{p-1}(\psi(x)-\psi(y))}{|x-y|^{n+ps}}dxdy
\end{align}
and the norm
\begin{align}\label{P3}
\|u\|_{X_{0}}=\left(\int_{Q}\frac{|u(x)-u(y)|^{p}}{|x-y|^{n+ps}}dxdy\right)^{\frac{1}{p}}
\end{align}
is equivalent to the usual one defined in \eqref{P1}. Since $u=0$ a.e. in $\mathbb{R}^{n}\setminus\Omega$, we have the \eqref{P2}, \eqref{P4} and \eqref{P3} can be extended to all $\mathbb{R}^{n}$. By results of \cite{SK,DPV}, the embedding $X_{0}\hookrightarrow L^{r}(\Omega)$ is continuous for any $r\in[1,p^{\star}]$ and compact whenever $r\in[1,p^{\star})$. Let $S$ be the best Sobolev constant for the embedding of $X_{0}\hookrightarrow L^{\alpha+\beta}(\Omega)$  defined by $$ S=\inf\limits_{z\in X_{0}\setminus{0}}\frac{\int_{Q}\frac{|u(x)-u(y)|^{p}}{|x-y|^{n+ps}}dxdy}{\left(\int_{\Omega}|u(x)|^{\alpha+\beta}dx\right)^{\frac{p}{\alpha+\beta}}}.$$
For further details on $X$ and $X_{0}$ and also for their properties we refer \cite{DPV} and the references therein. Let $E=X_{0}\times X_{0}$ be the Cartesian product of two Hilbert spaces, which is a reflexive Banach space endowed with the norm
\begin{align}\label{P5}
|(u,v)\|=\left(\|u\|^{p}_{X_{0}}+\|v\|^{p}_{X_{0}}\right)^{\frac{1}{p}}=\left(\int_{Q}\frac{|u(x)-u(y)|^{p}}{|x-y|^{n+ps}}dxdy+\int_{Q}\frac{|v(x)-v(y)|^{p}}{|x-y|^{n+ps}}dxdy\right)^{\frac{1}{p}}.
\end{align}

\begin{definition}
We say that $(u,v)\in E$ is a weak solution of \eqref{int1} if the identity
\begin{align*}
&\int_{Q}\frac{|u(x)-u(y)|^{p-2}(u(x)-u(y))(\varphi_{1}(x)-\varphi_{1}(y))}{|x-y|^{n+ps}}dxdy\\
&+\int_{Q}\frac{|v(x)-v(y)|^{p-2}(v(x)-v(y))(\varphi_{2}(x)-\varphi_{2}(y))}{|x-y|^{n+ps}}dxdy\\
&=\int_{\Omega}\left(\lambda f|u|^{q-2}u\varphi_{1}+\mu g|v|^{q-2}v\varphi_{2}\right)dx+\frac{2\alpha}{\alpha+\beta}\int_{\Omega}h|u|^{\alpha-2}u|v|^{\beta}\varphi_{1}dx+\frac{2\beta}{\alpha+\beta}\int_{\Omega}h|u|^{\alpha}|v|^{\beta-2}v\varphi_{2}dx.
\end{align*}
holds for all $(\varphi_{1},\varphi_{2})\in E$.
\end{definition}
Note that, the energy functional associated with \eqref{int1} is given by
\begin{align}
I_{\lambda,\mu}(u,v):&=\frac{1}{p}\int_{Q}\frac{|u(x)-u(y)|^{p}}{|x-y|^{n+ps}}dxdy+\frac{1}{p}\int_{Q}\frac{|v(x)-v(y)|^{p}}{|x-y|^{n+ps}}dxdy\\\nonumber
&-\frac{1}{q}\int_{\Omega}\left(\lambda f|u|^{q}+\mu g|v|^{q}\right)dx-\frac{2}{\alpha+\beta}\int_{\Omega}h|u|^{\alpha}|v|^{\beta}dx.
\end{align}

In the end of this section, we recall some notations that will be used in the sequel.\\
$\bullet$ $L^{p}(\Omega), \ 1\leq p\leq \infty$ denotes Lebesgue space with norm $\|\cdot\|_{p}$ and $E=X_{0}\times X_{0}$ is equipped with the norm $\|(u,v)\|=\left(\|u\|^{p}_{X_{0}}+\|v\|^{p}_{X_{0}}\right)^{\frac{1}{p}}=\left(\int_{Q}\frac{|u(x)-u(y)|^{p}}{|x-y|^{n+ps}}dxdy+\int_{Q}\frac{|v(x)-v(y)|^{p}}{|x-y|^{n+ps}}dxdy\right)^{\frac{1}{p}}.$\\
$\bullet$ The dual space of a Banach space $E$ will be denoted by $E^{-1}$. We set $t(u,v)=(tu,tv)$ for all $(u,v)\in E$ and $t\in \mathbb{R}$, $z=(u,v)$ is said to be positive if $u(x,y)>0,\ v(x,y)>0 $ in $E$ and to be non-negative if $u(x,y)\geq0,\ v(x,y)\geq0 $ in $E$.\\
$\bullet$ $B(0;r)$ is the ball at the origin with radius $r$. $o_{n}(1)$ denotes $o_{n}(1)\rightarrow0$ as $n\rightarrow+\infty$.\\
$\bullet$ $C, C_{i}, c$ will denote various positive constants which may vary from line to line.

\section{ The Nehari manifold }\label{sec3}
We consider the Nehari minimization problem: for $(\lambda,\mu)\in \mathbb{R}^{2}\setminus{(0,0)}$,
$$\theta_{\lambda,\mu}=\inf\{I_{\lambda,\mu}(u,v)\mid (u,v)\in\mathbf{\mathbb}{N}_{\lambda,\mu}\}$$
where $\mathbf{\mathbb}{N}_{\lambda,\mu}:=\{(u,v)\in E \setminus\{(0,0)\}\mid \langle I'_{\lambda,\mu}(u,v),(u,v)\rangle =0\}$ and
\begin{align}\label{N4}
\langle I'_{\lambda,\mu}(u,v),(u,v)\rangle=\|(u,v)\|^{p}-\int_{\Omega}\left(\lambda f|u|^{q}+\mu g|v|^{q}\right)dx-2\int_{\Omega}h|u|^{\alpha}|v|^{\beta}dx.
\end{align}
Note that $\mathbf{\mathbb}{N}_{\lambda,\mu}$ contains every nonzero solution of problem \eqref{int1}.

Define $$\langle \Phi_{\lambda,\mu}(u,v),(u,v)\rangle=\langle I'_{\lambda,\mu}(u,v),(u,v)\rangle.$$
Then $$\langle \Phi'_{\lambda,\mu}(u,v),(u,v)\rangle=p\|(u,v)\|^{p}-q\int_{\Omega}\left(\lambda f|u|^{q}+\mu g|v|^{q}\right)dx-2(\alpha+\beta)\int_{\Omega}h|u|^{\alpha}|v|^{\beta}dx.$$
Moreover, if $\int_{\Omega}\left(\lambda f|u|^{q}+\mu g|v|^{q}\right)dx\neq 0$ and $(u,v)\in \mathbf{\mathbb}{N}_{\lambda,\mu}$, we have
\begin{align}\label{N3}
\langle \Phi'_{\lambda,\mu}(u,v),(u,v)\rangle=(p-q)\|(u,v)\|^{p}-2(\alpha+\beta-q)\int_{\Omega}h|u|^{\alpha}|v|^{\beta}dx.
\end{align}
Similarly to the method used in \cite{TF}, we split $\mathbf{\mathbb}{N}_{\lambda,\mu}$ into three parts.
$$\mathbf{\mathbb}{N}^{+}_{\lambda,\mu}=\{(u,v)\in \mathbf{\mathbb}{N}_{\lambda,\mu}\mid \langle \Phi'_{\lambda,\mu}(u,v),(u,v)\rangle >0\};$$
$$\mathbf{\mathbb}{N}^{0}_{\lambda,\mu}=\{(u,v)\in \mathbf{\mathbb}{N}_{\lambda,\mu}\mid \langle \Phi'_{\lambda,\mu}(u,v),(u,v)\rangle =0\};$$
$$\mathbf{\mathbb}{N}^{-}_{\lambda,\mu}=\{(u,v)\in \mathbf{\mathbb}{N}_{\lambda,\mu}\mid \langle \Phi'_{\lambda,\mu}(u,v),(u,v)\rangle <0\}.$$
Then, we have the following results
\begin{lemma}\label{L1}
For each $(\lambda,\mu)\in \Theta$,
we have $\mathbf{\mathbb}{N}^{0}_{\lambda,\mu}=\emptyset$.
\end{lemma}
\begin{proof}
We consider the following two cases\\
{\bf Case 1}: $(u,v)\in \mathbf{\mathbb}{N}_{\lambda,\mu}$ and $\int_{\Omega}h|u|^{\alpha}|v|^{\beta}dx\leq 0$, we have$$\int_{\Omega}\left(\lambda f|u|^{q}+\mu g|v|^{q}\right)dx=\|(u,v)\|^{p}-2\int_{\Omega}h|u|^{\alpha}|v|^{\beta}dx>0.$$
Thus $\langle  \Phi'_{\lambda,\mu}(u,v),(u,v)\rangle=(p-q)\|(u,v)\|^{p}-2(\alpha+\beta-q)\int_{\Omega}h|u|^{\alpha}|v|^{\beta}dx>0$ \\ and so $(u,v)\notin \mathbf{\mathbb}{N}^{0}_{\lambda,\mu}$.\\
{\bf Case 2}: $(u,v)\in \mathbf{\mathbb}{N}_{\lambda,\mu}$ and $\int_{\Omega}h|u|^{\alpha}|v|^{\beta}dx> 0$.
Suppose that $\mathbf{\mathbb}{N}^{0}_{\lambda,\mu}\neq\emptyset$ for all $(\lambda,\mu)\in \mathbb{R}^{2}\setminus{(0,0)}$. Then for each $(u,v)\in \mathbf{\mathbb}{N}^{0}_{\lambda,\mu}$, we have
\begin{align}\label{L4}
\langle \Phi'_{\lambda,\mu}(u,v),(u,v)\rangle=(p-q)\|(u,v)\|^{p}-2(\alpha+\beta-q)\int_{\Omega}h|u|^{\alpha}|v|^{\beta}dx=0.
\end{align}
Thus
\begin{align*}
\|(u,v)\|^{p}=\frac{2(\alpha+\beta-q)}{p-q}\int_{\Omega}h|u|^{\alpha}|v|^{\beta}dx
\end{align*}
and
\begin{align*}
\int_{\Omega}\left(\lambda f|u|^{q}+\mu g|v|^{q}\right)dx&=\|(u,v)\|^{p}-2\int_{\Omega}h|u|^{\alpha}|v|^{\beta}dx\\\nonumber
&=\frac{2(\alpha+\beta-p)}{p-q}\int_{\Omega}h|u|^{\alpha}|v|^{\beta}dx>0.
\end{align*}
By the H\"{o}lder and Sobolev inequalities, we have
\begin{align}\label{N1}
\|(u,v)\|\geq [\frac{p-q}{2(\alpha+\beta-q)}S^{\frac{\alpha+\beta}{p}}]^{\frac{1}{\alpha+\beta-p}}
\end{align}
and
\begin{align*}
\frac{\alpha+\beta-p}{\alpha+\beta-q}\|(u,v)\|^{p}&=\|(u,v)\|^{p}-2\int_{\Omega}h|u|^{\alpha}|v|^{\beta}dx=\int_{\Omega}\left(\lambda f|u|^{q}+\mu g|v|^{q}\right)dx\\
&\leq |\lambda|\|f\|_{L^{q^{\star}}}|\|u\|^{q}_{L^{\alpha+\beta}}+|\mu|\|g\|_{L^{q^{\star}}}|\|v\|^{q}_{L^{\alpha+\beta}}\\
&\leq [(|\lambda|\|f\|_{L^{q^{\star}}})^{\frac{p}{p-q}}+(|\mu|\|g\|_{L^{q^{\star}}})^{\frac{p}{p-q}}]^{\frac{p-q}{p}}S^{-\frac{q}{p}}\|(u,v)\|^{q}.
\end{align*}
This implies
\begin{align}\label{N2}
\|(u,v)\|\leq \left(S^{-\frac{q}{p}}\frac{\alpha+\beta-q}{\alpha+\beta-p}\right)^{\frac{1}{p-q}}\left[(|\lambda|\|f\|_{L^{q^{\star}}})^{\frac{p}{p-q}}+(|\mu|\|g\|_{L^{q^{\star}}})^{\frac{p}{p-q}}\right]^{\frac{1}{p}}.
\end{align}
By \eqref{N1} and \eqref{N2}, we have
$$[(|\lambda|\|f\|_{L^{q^{\star}}})^{\frac{p}{p-q}}+(|\mu|\|g\|_{L^{q^{\star}}})^{\frac{p}{p-q}}]\geq\left[\frac{p-q}{2(\alpha+\beta-p)}S^{\frac{\alpha+\beta}{p}}\right]^{\frac{p}{\alpha+\beta-p}}\left(S^{-\frac{q}{p}}\frac{\alpha+\beta-q}{\alpha+\beta-p}\right)^{-\frac{p}{p-q}},$$
contradicting with the assumption.
\end{proof}
Lemma \ref{L1} suggests that for each $(\lambda,\mu)\in \Theta$, we can write $\mathbf{\mathbb}{N}_{\lambda,\mu}=\mathbf{\mathbb}{N}^{+}_{\lambda,\mu}\cup \mathbf{\mathbb}{N}^{-}_{\lambda,\mu}$.

 Next, we define $$\theta^{+}_{\lambda,\mu}=\inf\limits_{z\in \mathbf{\mathbb}{N}^{+}_{\lambda,\mu}}I_{\lambda,\mu}(z)\ \ and \ \  \theta^{-}_{\lambda,\mu}=\inf\limits_{z\in \mathbf{\mathbb}{N}^{-}_{\lambda,\mu}}I_{\lambda,\mu}(z).$$
The following Lemma shows that the minimizer on $\mathbf{\mathbb}{N}_{\lambda,\mu}$ is critical point for $I_{\lambda,\mu}$
\begin{lemma}\label{L22}
For each $(\lambda,\mu)\in\Theta$, let $(u_{0},v_{0})$ be a local minimizer for $I_{\lambda,\mu}$ on $\mathbf{\mathbb}{N}_{\lambda,\mu}$ ,then $I'_{\lambda,\mu}(u_{0},v_{0})=0$ in $E^{-1}$.
\begin{proof}
 Since $(u_{0},v_{0})$ is a local minimizer for $I_{\lambda,\mu}$ on $\mathbf{\mathbb}{N}_{\lambda,\mu}$, that is $(u_{0},v_{0})$ is a solution of the optimization problem$$\min\{I_{\lambda,\mu}(u,v)\mid \Phi_{\lambda,\mu}(u,v)=0 \}.$$
Then, by the theory of Lagrange multipliers, there exists a constant $L\in\mathbb{R}$ such that$$\langle I'_{\lambda,\mu}(u_{0},v_{0}),(u_{0},v_{0})\rangle=L\langle\Phi'_{\lambda,\mu}(u_{0},v_{0}),(u_{0},v_{0})\rangle.$$
Since $(u_{0},v_{0})\notin \mathbf{\mathbb}{N}^{0}_{\lambda,\mu}$, we have $\langle \Phi'_{\lambda,\mu}(u_{0},v_{0}),(u_{0},v_{0})\rangle \neq0$, thus $L=0$, this completes the proof.
\end{proof}
\end{lemma}
Moreover, we have the following properties about the Nehari manifold $\mathbf{\mathbb}{N}_{\lambda,\mu}.$
\begin{lemma}\label{L2}
we have\\
$(i)$ If $(u,v)\in \mathbf{\mathbb}{N}^{+}_{\lambda,\mu}$, then $\int_{\Omega}\left(\lambda f|u|^{q}+\mu g|v|^{q}\right)dx>0$\\
$(ii)$ If $(u,v)\in \mathbf{\mathbb}{N}^{-}_{\lambda,\mu}$, then $\int_{\Omega}h|u|^{\alpha}|v|^{\beta}dx>0$
\end{lemma}
\begin{proof}
$(i)$ We consider the following two cases.\\
{\bf Case 1}: If $\int_{\Omega}h|u|^{\alpha}|v|^{\beta}dx\leq0$, we have
$$\int_{\Omega}\left(\lambda f|u|^{q}+\mu g|v|^{q}\right)dx=\|(u,v)\|^{p}-2\int_{\Omega}h|u|^{\alpha}|v|^{\beta}dx>0.$$
{\bf Case 2}: If $\int_{\Omega}h|u|^{\alpha}|v|^{\beta}dx>0$, since
 $$\|(u,v)\|^{p}-\int_{\Omega}\left(\lambda f|u|^{q}+\mu g|v|^{q}\right)dx-2\int_{\Omega}h|u|^{\alpha}|v|^{\beta}dx=0$$
 and
 $$\langle \Phi'_{\lambda,\mu}(u,v),(u,v)\rangle=p\|(u,v)\|^{p}-q\int_{\Omega}\left(\lambda f|u|^{q}+\mu g|v|^{q}\right)dx-2(\alpha+\beta)\int_{\Omega}h|u|^{\alpha}|v|^{\beta}dx>0,$$
 it follows that $$(p-q)\int_{\Omega}\left(\lambda f|u|^{q}+\mu g|v|^{q}\right)dx-2(\alpha+\beta-q)\int_{\Omega}h|u|^{\alpha}|v|^{\beta}dx>0,$$
 which implies
 $$\int_{\Omega}\left(\lambda f|u|^{q}+\mu g|v|^{q}\right)dx>\frac{2(\alpha+\beta-q)}{p-q}\int_{\Omega}h|u|^{\alpha}|v|^{\beta}dx>0.$$
 $(ii)$ We consider the following two cases.\\
 {\bf Case 1}: If $\int_{\Omega}\left(\lambda f|u|^{q}+\mu g|v|^{q}\right)dx=0$, we have
 $$2\int_{\Omega}h|u|^{\alpha}|v|^{\beta}dx=\|(u,v)\|^{p}>0.$$
  {\bf Case 2}: If $\int_{\Omega}\left(\lambda f|u|^{q}+\mu g|v|^{q}\right)dx\neq0$, we have
  $$(p-q)\|(u,v)\|^{p}-2(\alpha+\beta-q)\int_{\Omega}h|u|^{\alpha}|v|^{\beta}dx=\langle \Phi'_{\lambda,\mu}(u,v),(u,v)\rangle<0.$$
  Thus $\int_{\Omega}h|u|^{\alpha}|v|^{\beta}dx>0$.
\end{proof}
\begin{lemma}\label{hhh}
The following facts hold\\
$(i)$ If $(\lambda,\mu)\in\Theta$, then we have $\theta_{\lambda,\mu}\leq\theta^{+}_{\lambda,\mu}<0$\\
$(ii$ If $(\lambda,\mu)\in\Psi$, then we have $\theta^{-}_{\lambda,\mu}>c_{0}$ for some positive constant $c_{0}$ depending on $\lambda,\mu,p,q,S,$.\\
$(iii)$ The energy functional $I_{\lambda,\mu}$ is bounded below and coercive on $\mathbf{\mathbb}{N}_{\lambda,\mu}$.
\end{lemma}
\begin{proof}
$(i)$ Let $(u,v)\in \mathbf{\mathbb}{N}^{+}_{\lambda,\mu}$, by \eqref{N3}, we have
$$\frac{p-q}{2(\alpha+\beta-q)}\|(u,v)\|^{p}>\int_{\Omega}h|u|^{\alpha}|v|^{\beta}dx.$$
Hence
\begin{align*}
I_{\lambda,\mu}(u,v)=&\left(\frac{1}{p}-\frac{1}{q}\right)\|(u,v)\|^{p}+2\left(\frac{1}{q}-\frac{1}{\alpha+\beta}\right)\int_{\Omega}h|u|^{\alpha}|v|^{\beta}dx\\
&\leq\left[\left(\frac{1}{p}-\frac{1}{q}\right)+\left(\frac{1}{q}-\frac{1}{\alpha+\beta}\right)\frac{p-q}{\alpha+\beta-q}\right]\|(u,v)\|^{p}\\
&\leq\frac{(q-p)(\alpha+\beta-p)}{pq(\alpha+\beta)}\|(u,v)\|^{p}<0.
\end{align*}
Therefore, by the definition of $\theta_{\lambda,\mu},\ \theta^{+}_{\lambda,\mu}$, we can deduce that $\theta_{\lambda,\mu}\leq\theta^{+}_{\lambda,\mu}<0$.\\
$(ii)$ Let $(u,v)\in \mathbf{\mathbb}{N}^{-}_{\lambda,\mu}$, by \eqref{N3}, we have
$$\frac{p-q}{2(\alpha+\beta-q)}\|(u,v)\|^{p}<\int_{\Omega}h|u|^{\alpha}|v|^{\beta}dx.$$
By the H\"{o}lder inequality and the Sobolev embedding theorem, we have
$$\int_{\Omega}h|u|^{\alpha}|v|^{\beta}dx\leq S^{-\frac{\alpha+\beta}{p}}\|(u,v)\|^{\alpha+\beta}.$$
Hence
\begin{align}\label{TH4}
\|(u,v)\|>\left(\frac{p-q}{2(\alpha+\beta-q)}S^{\frac{\alpha+\beta}{p}}\right)^{\frac{1}{\alpha+\beta-p}}\ for\  all\ z\in\mathbf{\mathbb}{N}^{-}_{\lambda,\mu}.
\end{align}
By \eqref{TH4}, we have
\begin{align*}
I_{\lambda,\mu}(u,v)&=\frac{\alpha+\beta-p}{p(\alpha+\beta)}\|(u,v)\|^{p}-\frac{\alpha+\beta-q}{q(\alpha+\beta)}\int_{\Omega}\left(\lambda f|u|^{q}+\mu g|v|^{q}\right)dx\\
&\geq \|(u,v)\|^{q}\left[\frac{\alpha+\beta-p}{p(\alpha+\beta)}\|(u,v)\|^{p-q}-\frac{\alpha+\beta-q}{q(\alpha+\beta)}[(|\lambda|\|f\|_{L^{q^{\star}}})^{\frac{p}{p-q}}+(|\mu|\|g\|_{L^{q^{\star}}})^{\frac{p}{p-q}}]^{\frac{p-q}{p}}S^{-\frac{q}{p}}\right]\\
&>\left\{-\frac{\alpha+\beta-q}{q(\alpha+\beta)}[(|\lambda|\|f\|_{L^{q^{\star}}})^{\frac{p}{p-q}}+(|\mu|\|g\|_{L^{q^{\star}}})^{\frac{p}{p-q}}]^{\frac{p-q}{p}}S^{-\frac{q}{p}} +\right.\\
&\phantom{=\;\;}\left.\frac{\alpha+\beta-p}{p(\alpha+\beta)}\left(\frac{p-q}{2(\alpha+\beta-q)}S^{\frac{\alpha+\beta}{p}}\right)^{\frac{p-q}{\alpha+\beta-p}} \right\} \times \left(\frac{p-q}{2(\alpha+\beta-q)}S^{\frac{\alpha+\beta}{p}}\right)^{\frac{q}{\alpha+\beta-p}}.
\end{align*}
Thus, if $(\lambda,\mu)\in\Psi$, then $$I_{\lambda,\mu}>c_{0}, \ for \ all \ z\in\mathbf{\mathbb}{N}^{-}_{\lambda,\mu},$$
for some positive constant $c_{0}=c_{0}(\lambda,\mu,p,q,S)$.\\
$(iii)$ Let $(u,v)\in \mathbf{\mathbb}{N}_{\lambda,\mu}$, by \eqref{N4}, H\"{o}lder and Sobolev inequality, we have
\begin{align*}
I_{\lambda,\mu}(u,v)&=\frac{\alpha+\beta-p}{p(\alpha+\beta)}\|(u,v)\|^{p}-\frac{\alpha+\beta-q}{q(\alpha+\beta)}\int_{\Omega}\left(\lambda f|u|^{q}+\mu g|v|^{q}\right)dx\\
&\geq \frac{\alpha+\beta-p}{p(\alpha+\beta)}\|(u,v)\|^{p}-\frac{\alpha+\beta-q}{q(\alpha+\beta)}[(|\lambda|\|f\|_{L^{q^{\star}}})^{\frac{p}{p-q}}+(|\mu|\|g\|_{L^{q^{\star}}})^{\frac{p}{p-q}}]^{\frac{p-q}{p}}S^{-\frac{q}{p}} \|(u,v)\|^{q}.
\end{align*}
Since $1<q<p$, then the energy functional $I_{\lambda,\mu}$ is bounded below and coercive on $\mathbf{\mathbb}{N}_{\lambda,\mu}$.
\end{proof}
For each $(\lambda,\mu)\in \mathbf{\mathbb}{N}^{-}_{\lambda,\mu}$, we write
$$t_{max}=\left(\frac{(p-q)\|(u,v)\|^{p}}{2(\alpha+\beta-q)\int_{\Omega}h|u|^{\alpha}|v|^{\beta}dx}\right)^{\frac{1}{\alpha+\beta-p}}>0.$$
Then the following Lemma holds.
\begin{lemma}\label{L21}
For each $(\lambda,\mu)\in\Theta$ and $(u,v)\in \mathbf{\mathbb}{N}^{-}_{\lambda,\mu}$, we have \\
$(i)$ If $\int_{\Omega}\left(\lambda f|u|^{q}+\mu g|v|^{q}\right)dx\leq0$, then there exists a unique  $(t^{-}u,t^{-}v)>0$ such that $(t^{-}u,t^{-}v)\in \mathbf{\mathbb}{N}^{-}_{\lambda,\mu}$ and $I_{\lambda,\mu}(t^{-}u,t^{-}v)=\max \limits_{t>0}I_{\lambda,\mu}(tu,tv)$.\\
$(ii)$ If $\int_{\Omega}\left(\lambda f|u|^{q}+\mu g|v|^{q}\right)dx>0$, then there exist unique $0<t^{+}<t_{max}<t^{-}$, such that $(t^{+}u,t^{+}v)\in \mathbf{\mathbb}{N}^{+}_{\lambda,\mu},\ (t^{-}u,t^{-}v)\in \mathbf{\mathbb}{N}^{-}_{\lambda,\mu}$ and $I_{\lambda,\mu}(t^{+}u,t^{+}v)=\min \limits_{0<t<t_{max}}I_{\lambda,\mu}(tu,tv),$\\
 $I_{\lambda,\mu}(t^{-}u,t^{-}v)=\max\limits_{t\geq 0}I_{\lambda,\mu}(tu,tv)$.
\end{lemma}
\begin{proof}
Fix $(u,v)\in \mathbf{\mathbb}{N}^{-}_{\lambda,\mu}$, by Lemma \ref{L2}, we have $\int_{\Omega}h|u|^{\alpha}|v|^{\beta}dx>0$. Let
$$m(t)=t^{p-q}\|(u,v)\|^{p}-2t^{\alpha+\beta-q}\int_{\Omega}h|u|^{\alpha}|v|^{\beta}dx,\  for \ t\geq0.$$
Clearly, $m(0)=0$, $m(t)\rightarrow-\infty$ as $t\rightarrow\infty$. Since $$m'(t)=(p-q)t^{p-q-1}\|(u,v)\|^{p}-2(\alpha+\beta-q)t^{\alpha+\beta-q-1}\int_{\Omega}h|u|^{\alpha}|v|^{\beta}dx,$$
we have that $m(t)$ is increasing for $t\in[0,t_{max})$, decreasing for $t\in(t_{max},+\infty)$ and achieves its maximum at $t_{max}$. Moreover,
\begin{align*}
m(t_{max})&=\left(\frac{(p-q)\|(u,v)\|^{p}}{2(\alpha+\beta-q)\int_{\Omega}h|u|^{\alpha}|v|^{\beta}dx}\right)^{\frac{p-q}{\alpha+\beta-p}}\|(u,v)\|^{p}\\
&-2\left(\frac{(p-q)\|(u,v)\|^{p}}{2(\alpha+\beta-q)\int_{\Omega}h|u|^{\alpha}|v|^{\beta}dx}\right)^{\frac{\alpha+\beta-q}{\alpha+\beta-p}}\int_{\Omega}h|u|^{\alpha}|v|^{\beta}dx\\
&=\|(u,v)\|^{q}\left[\left(\frac{p-q}{2(\alpha+\beta-q)}\right)^{\frac{p-q}{\alpha+\beta-p}}-2\left(\frac{p-q}{2(\alpha+\beta-q)}\right)^{\frac{\alpha+\beta-q}{\alpha+\beta-p}}\right]\left(\frac{\|(u,v)\|^{\alpha+\beta}}{\int_{\Omega}h|u|^{\alpha}|v|^{\beta}dx}\right)^{\frac{p-q}{\alpha+\beta-p}}\\
&\geq \|(u,v)\|^{q}\left(\frac{\alpha+\beta-p}{\alpha+\beta-q}\right)\left(S^{\frac{\alpha+\beta}{p}}\frac{p-q}{2(\alpha+\beta-q)}\right)^{\frac{p-q}{\alpha+\beta-p}}.
\end{align*}
That is
\begin{align}\label{L3}
m(t_{max})\geq \|(u,v)\|^{q}\left(\frac{\alpha+\beta-p}{\alpha+\beta-q}\right)\left(S^{\frac{\alpha+\beta}{p}}\frac{p-q}{2(\alpha+\beta-q)}\right)^{\frac{p-q}{\alpha+\beta-p}}.
\end{align}
$(i)$ If $\int_{\Omega}\left(\lambda f|u|^{q}+\mu g|v|^{q}\right)dx\leq0$, by the property of $m(t)$, there is a unique $t^{-}>t_{max}$ such that $m(t^{-})=\int_{\Omega}\left(\lambda f|u|^{q}+\mu g|v|^{q}\right)dx$ and $m'(t^{-})<0$.\\
Since
\begin{align}\label{L5}
&\langle \Phi'_{\lambda,\mu}(t^{-}u,t^{-}v),(t^{-}u,t^{-}v)\rangle=(p-q)(t^{-})^{p}\|(u,v)\|^{p}-2(\alpha+\beta-q)(t^{-})^{\alpha+\beta}\int_{\Omega}h|u|^{\alpha}|v|^{\beta}dx\\\nonumber
&=(t^{-})^{1+q}\left[(p-q)(t^{-})^{p-q-1}\|(u,v)\|^{p}-2(\alpha+\beta-q)(t^{-})^{\alpha+\beta-q-1}\int_{\Omega}h|u|^{\alpha}|v|^{\beta}dx\right]\\\nonumber
&=(t^{-})^{1+q}m'(t^{-})<0
\end{align}
and
\begin{align*}
\langle I'_{\lambda,\mu}(t^{-}u,t^{-}v),(t^{-}u&,t^{-}v)\rangle=(t^{-})^{p}\|(u,v)\|^{p}-(t^{-})^{q}\int_{\Omega}\left(\lambda f|u|^{q}+\mu g|v|^{q}\right)dx-2(t^{-})^{\alpha+\beta}\int_{\Omega}h|u|^{\alpha}|v|^{\beta}dx\\
&=(t^{-})^{q}\left[(t^{-})^{p-q}\|(u,v)\|^{p}-\int_{\Omega}\left(\lambda f|u|^{q}+\mu g|v|^{q}\right)dx-2(t^{-})^{\alpha+\beta-q}\int_{\Omega}h|u|^{\alpha}|v|^{\beta}dx\right]\\
&=(t^{-})^{q}\left[m(t^{-})-\int_{\Omega}\left(\lambda f|u|^{q}+\mu g|v|^{q}\right)dx\right]=0.
\end{align*}
Thus $(t^{-}u,t^{-}v)\in \mathbf{\mathbb}{N}^{-}_{\lambda,\mu}$.\\
For $t>t_{max}$, by \eqref{L5}, we know
$$(p-q)t^{p}\|(u,v)\|^{p}-2(\alpha+\beta-q)t^{\alpha+\beta}\int_{\Omega}h|u|^{\alpha}|v|^{\beta}dx<0.$$
When $tz\in \mathbf{\mathbb}{N}_{\lambda,\mu}$, we have
$$\|(u,v)\|^{p}-t^{q-p}\int_{\Omega}\left(\lambda f|u|^{q}+\mu g|v|^{q}\right)dx-2t^{\alpha+\beta-p}\int_{\Omega}h|u|^{\alpha}|v|^{\beta}dx=0$$
and
$$\frac{d^{2}}{dt^{2}}I_{\lambda,\mu}(tu,tv)=(p-1)t^{p-2}\|(u,v)\|^{p}-(q-1)t^{q-2}\int_{\Omega}\left(\lambda f|u|^{q}+\mu g|v|^{q}\right)dx-2(\alpha+\beta-1)t^{\alpha+\beta-2}\int_{\Omega}h|u|^{\alpha}|v|^{\beta}dx.$$
Consequently, $\frac{d^{2}}{dt^{2}}I_{\lambda,\mu}(tu,tv)=t^{q-1}m'(t)<0$.\\
Since$$\frac{d}{dt}I_{\lambda,\mu}(tu,tv)=t^{p-1}\|(u,v)\|^{p}-t^{q-1}\int_{\Omega}\left(\lambda f|u|^{q}+\mu g|v|^{q}\right)dx-2t^{\alpha+\beta-1}\int_{\Omega}h|u|^{\alpha}|v|^{\beta}dx.$$
We have $\frac{d}{dt}I_{\lambda,\mu}(tu,tv)=0$ for $t=t^{-}$.
Thus, $I_{\lambda,\mu}(t^{-}z)=\max \limits_{t>0}I_{\lambda,\mu}(tz)$.\\
$(ii)$ If $\int_{\Omega}\left(\lambda f|u|^{q}+\mu g|v|^{q}\right)dx>0$.\\
Since
\begin{align*}
m(0)&=0<\int_{\Omega}\left(\lambda f|u|^{q}+\mu g|v|^{q}\right)dx\\
&\leq [(|\lambda|\|f\|_{L^{q^{\star}}})^{\frac{p}{p-q}}+(|\mu|\|g\|_{L^{q^{\star}}})^{\frac{p}{p-q}}]^{\frac{p-q}{p}}S^{-\frac{q}{p}}\|(u,v)\|^{q}\\
&\leq \|(u,v)\|^{q}\left(\frac{\alpha+\beta-p}{\alpha+\beta-q}\right)\left(S^{\frac{\alpha+\beta}{p}}\frac{p-q}{2(\alpha+\beta-q)}\right)^{\frac{p-q}{\alpha+\beta-p}}\\
&\leq m(t_{max}) \ for \ (\lambda,\mu)\in\Theta,
\end{align*}
therefore, there are unique $t^{+}$ and $t^{-}$ such that $0<t^{+}<t_{max}<t^{-}$,
 $$m(t^{+})=\int_{\Omega}\left(\lambda f|u|^{q}+\mu g|v|^{q}\right)dx=m(t^{-})$$
 and $$m'(t^{+})>0>m'(t^{-}).$$
Thus, by the same arguments as $(i)$, we have  $(t^{+}u,t^{+}v)\in \mathbf{\mathbb}{N}^{+}_{\lambda,\mu}$,\ $(t^{-}u,t^{-}v)\in \mathbf{\mathbb}{N}^{-}_{\lambda,\mu}$,\ \ $I_{\lambda,\mu}(t^{-}u,t^{-}v)\geq I_{\lambda,\mu}(tu,tv)\geq I_{\lambda,\mu}(t^{+}u,t^{+}v)$ for each $t\in[t^{+},t^{-}]$ and $I_{\lambda,\mu}(t^{+}u,t^{+}v)\leq I_{\lambda,\mu}(tu,tv)$ for each $t\in[0,t^{+}]$.\\
That is
$$I_{\lambda,\mu}(t^{+}u,t^{+}v)=\min \limits_{0<t<t_{max}}I_{\lambda,\mu}(tu,tv),\ \ I_{\lambda,\mu}(t^{-}u,t^{-}v)=\max\limits_{t\geq 0}I_{\lambda,\mu}(tu,tv).$$
\end{proof}

\section{ Existence of Palais-Smale sequence }\label{sec4}
\begin{definition}
We say that $(u_{n},v_{n})\in E$ is a $(PS)_{c}$ sequence in $E$ for $I_{\lambda\mu}$, if $I_{\lambda\mu}(u_{n},v_{n})=c+o_{n}(1)$ and $I'_{\lambda\mu}(u_{n},v_{n})=o_{n}(1)$  strongly in $E^{-1}$ as $n\rightarrow\infty$. If any $(PS)_{c}$ sequence in $E$ for $I_{\lambda\mu}$ admits a convergent subsequence, we say that $I_{\lambda\mu}$ satisfies the $(PS)_{c}$ condition.
\end{definition}
First, we will use the idea of Tarantello \cite{T} to get the following results.
\begin{lemma}\label{L7}
Let $z=(u,v)$ and $(\lambda,\mu)\in\Theta$, then for each $z \in \mathbf{\mathbb}{N}_{\lambda,\mu}$, there exists $r>0$ and a differentiable function $\xi:B(0,r)\subset E\rightarrow \mathbb{R}^{+}$ such that $\xi(0)=1$ and $\xi(m)(z-m)\in \mathbf{\mathbb}{N}_{\lambda,\mu}$ for every $m\in B(0,r)$. Let
\begin{align*}
&\mathrm{T}_{1}:=p\int_{Q}\frac{|u(x)-u(y)|^{p-2}(u(x)-u(y))(m_{1}(x)-m_{1}(y))}{|x-y|^{n+ps}}dxdy\\
&+p\int_{Q}\frac{|v(x)-v(y)|^{p-2}(v(x)-v(y))(m_{2}(x)-m_{2}(y))}{|x-y|^{n+ps}}dxdy,\\
&\mathrm{T}_{2}:=q\int_{\Omega}\left(\lambda f|u|^{q-2}um_{1}+\mu g|v|^{q-2}vm_{2}\right)dx,\\
&\mathrm{T}_{3}:=2\int_{\Omega}\left(\alpha h|u|^{\alpha-2}um_{1}|v|^{\beta}+\beta h|u|^{\alpha}|v|^{\beta-2}vm_{2}\right)dx,
\end{align*}
then
\begin{align}\label{L6}
\langle \xi'(0),m\rangle=\frac{\mathrm{T}_{2}+\mathrm{T}_{3}-\mathrm{T}_{1}}{(p-q)\|(u,v)\|^{p}-2(\alpha+\beta-q)\int_{\Omega}h|u|^{\alpha}|v|^{\beta}dx}
\end{align}
holds for all $m\in E.$
\end{lemma}
\begin{proof}
For $z\in \mathbf{\mathbb}{N}_{\lambda,\mu}$, define a function $F:\mathbb{R}\times E\rightarrow\mathbb{R}$ by
\begin{align*}
&F_{z}(\xi,q):=\langle I'_{\lambda,\mu}(\xi(z-q)),\xi(z-q)\rangle\\
&=\xi^{p}\int_{Q}\frac{|(u(x)-q_{1}(x))-(u(y)-q_{1}(y))|^{p}}{|x-y|^{n+ps}}dxdy+\xi^{p}\int_{Q}\frac{|(v(x)-q_{2}(x))-(v(y)-q_{2}(y))|^{p}}{|x-y|^{n+ps}}dxdy\\
&-\xi^{q}\int_{\Omega}\left(\lambda f|u-q_{1}|^{q}+\mu g|v-q_{2}|^{q}\right)dx-2\xi^{\alpha+\beta}\int_{\Omega}h|u-q_{1}|^{\alpha}|v-q_{2}|^{\beta}dx.
\end{align*}
Then, $F_{z}(1,0)=\langle I'_{\lambda,\mu}(u,v),(u,v)\rangle=0$ and by Lemma \ref{L1}, we have  $\mathbf{\mathbb}{N}^{0}_{\lambda,\mu}=\emptyset$.\\ That is
\begin{align*}
\frac{dF_{z}(1,0)}{d\xi}=&p\|(u,v)\|^{p}-q\int_{\Omega}\left(\lambda f|u|^{q}+\mu g|v|^{q}\right)dx-2(\alpha+\beta)\int_{\Omega}h|u|^{\alpha}|v|^{\beta}dx\\
&=(p-q)\|(u,v)\|^{p}-2(\alpha+\beta-q)\int_{\Omega}h|u|^{\alpha}|v|^{\beta}dx\neq0.
\end{align*}
According to the implicit function theorem, there exist $r>0$ and a differentiable function $\xi:B(0,r)\subset E\rightarrow \mathbb{R}^{+}$ such that $\xi(0)=1$ and \eqref{L6} holds. Moreover, $F_{z}(\xi(m),m)=0$ holds for all $m\in B(0,r)$ is equivalent to $\langle I'_{\lambda,\mu}(\xi(m)(z-m)),\xi(m)(z-m)\rangle=0$ for all $m\in B(0,r)$.
That is $\xi(m)(z-m)\in \mathbf{\mathbb}{N}_{\lambda,\mu}$.
\end{proof}
\begin{lemma}\label{M3}
Let $z=(u,v)$ and $(\lambda,\mu)\in\Theta$, then for each $z\in \mathbf{\mathbb}{N}^{-}_{\lambda,\mu}$, there exists $r>0$ and a differentiable function $\xi^{-}:B(0,r)\subset E\rightarrow \mathbb{R}^{+}$ such that $\xi^{-}(0)=1$ and $\xi^{-}(m)(z-m)\in \mathbf{\mathbb}{N}_{\lambda,\mu}$ for every $m\in B(0,r)$ and formula \eqref{L6} holds.
\end{lemma}
\begin{proof}
Similar to the argument in Lemma \ref{L7}, there exists $r>0$ and a differentiable function $\xi^{-}:B(0,r)\subset E\rightarrow \mathbb{R}^{+}$ such that $\xi^{-}(0)=1$ and $\xi^{-}(m)(z-m)\in \mathbf{\mathbb}{N}_{\lambda,\mu}$ for every $m\in B(0,r)$ and formula \eqref{L6} holds.
 Since
 $$\langle \Phi'_{\lambda,\mu}(z),z\rangle=(p-q)\|(u,v)\|^{p}-2(\alpha+\beta-q)\int_{\Omega}h|u|^{\alpha}|v|^{\beta}dx<0,$$
 by the continuity of function $\Phi'_{\lambda,\mu}$ and $\xi^{-}$, we have
 \begin{align*}
 &\langle \Phi'_{\lambda,\mu}(\xi^{-}(m)(z-m)),\xi^{-}(m)(z-m)\rangle\\
 &=(p-q)\|\xi^{-}(m)(z-m)\|^{2}-2(\alpha+\beta-q)\int_{\Omega}h|(\xi^{-}(m)(z-m))_{1}|^{\alpha}|(\xi^{-}(m)(z-m))_{2}|^{\beta}dx<0.
 \end{align*}
 This implies that $\xi^{-}(m)(z-m)\in \mathbf{\mathbb}{N}^{-}_{\lambda,\mu}.$
\end{proof}
\begin{lemma}\label{L10}
The following facts hold:\\
$(i)$ If $(\lambda,\mu)\in \Theta$, then there is a $(PS)_{\theta_{\lambda,\mu}}$-sequence $\{z_{n}\}=\{(u_{n},v_{n})\}\subset\mathbf{\mathbb}{N}_{\lambda,\mu}$ for $I_{\lambda,\mu}$.\\
$(ii)$  If $(\lambda,\mu)\in \Psi$, then there is a $(PS)_{\theta^{-}_{\lambda,\mu}}$-sequence $\{z_{n}\}=\{(u_{n},v_{n})\}\subset\mathbf{\mathbb}{N}^{-}_{\lambda,\mu}$ for $I_{\lambda,\mu}$.
\end{lemma}
\begin{proof}
(1) By Lemma \ref{hhh}(iii) and Ekeland Variational Principle \cite{IE}, there exists a minimizing sequence $\{z_{n}\}\subset \mathbf{\mathbb}{N}_{\lambda,\mu}$ such that
\begin{align}\label{defg}
&I_{\lambda,\mu}(z_{n})<\theta_{\lambda,\mu}+\frac{1}{n},\\\nonumber
&I_{\lambda,\mu}(z_{n})<I_{\lambda,\mu}(w)+\frac{1}{n}\|w-z_{n}\|, \ \forall \ w\in \mathbf{\mathbb}{N}_{\lambda,\mu}.
\end{align}
By taking $n$ large, from  Lemma \ref{hhh}(i), we have $\theta_{\lambda,\mu}<0$, thus
\begin{align}\label{ggg}
I_{\lambda,\mu}(z_{n})=&(\frac{1}{p}-\frac{1}{\alpha+\beta})\|(u_{n},v_{n})\|^{p}-(\frac{1}{q}-\frac{1}{\alpha+\beta})\int_{\Omega}\left(\lambda f|u_{n}|^{q}+\mu g|v_{n}|^{q}\right)dx\\\nonumber
&<\theta_{\lambda,\mu}+\frac{1}{n}<\frac{\theta_{\lambda,\mu}}{2}.
\end{align}
This implies
\begin{align}\label{iii}
-\frac{q(\alpha+\beta)}{2(\alpha+\beta-q)}\theta_{\lambda,\mu}&<\int_{\Omega}\left(\lambda f|u_{n}|^{q}+\mu g|v_{n}|^{q}\right)dx\\\nonumber
&\leq [(|\lambda|\|f\|_{L^{q^{\star}}})^{\frac{p}{p-q}}+(|\mu|\|g\|_{L^{q^{\star}}})^{\frac{p}{p-q}}]^{\frac{p-q}{p}}S^{-\frac{q}{p}}\|(u_{n},v_{n})\|^{q}.
\end{align}
Consequently, $z_{n}\neq0$ and putting together \eqref{ggg}, \eqref{iii} and the H\"{o}lder inequality, we obtain
\begin{align*}
\|(u_{n},v_{n})\|> \left[-\frac{q(\alpha+\beta)}{2(\alpha+\beta-q)}\theta_{\lambda,\mu}\left[(|\lambda|\|f\|_{L^{q^{\star}}})^{\frac{p}{p-q}}+(|\mu|\|g\|_{L^{q^{\star}}})^{\frac{p}{p-q}}\right]^{\frac{q-p}{p}}S^{\frac{q}{p}}\right]^{\frac{1}{q}}.
\end{align*}
and
\begin{align}\label{abc}
\|(u_{n},v_{n})\|< \left[\frac{p(\alpha+\beta-q)}{q(\alpha+\beta-2)}\left[(|\lambda|\|f\|_{L^{q^{\star}}})^{\frac{p}{p-q}}+(|\mu|\|g\|_{L^{q^{\star}}})^{\frac{p}{p-q}}\right]^{\frac{p-q}{p}}S^{-\frac{q}{p}}\right]^{\frac{1}{p-q}}.
\end{align}
Now, we will show that
$$\|I'_{\lambda,\mu}(z_{n})\|_{E^{-1}}\rightarrow0\ \  as\ \ n\rightarrow+\infty.$$
Applying Lemma \ref{L7} to $z_{n}$, we can obtain the function $\xi_{n}:B(0,r_{n})\subset E\rightarrow \mathbb{R}^{+}$ such that $\xi_{n}(0)=1$ and $\xi_{n}(m)(z_{n}-m)\in \mathbf{\mathbb}{N}_{\lambda,\mu}$ for every $m\in B(0,r_{n})$. Taking $0<\rho<r_{n}$, let $w\in E$ with $w\neq0$ and put $m^{\star}=\frac{\rho w}{\|w\|}$. We set $m_{\rho}=\xi_{n}(m^{\star})(z_{n}-m^{\star})$, then $m_{\rho}\in \mathbf{\mathbb}{N}_{\lambda,\mu}$. \\By \eqref{defg}, we have
\begin{align*}
I_{\lambda,\mu}(m_{\rho})-I_{\lambda,\mu}(z_{n})\geq -\frac{1}{n}\|m_{\rho}-z_{n}\|.
\end{align*}
By the Mean Value Theorem, we get
$$\langle I'_{\lambda,\mu}(z_{n}),m_{\rho}-z_{n}\rangle+o(\|m_{\rho}-z_{n}\|)\geq -\frac{1}{n}\|m_{\rho}-z_{n}\|.$$
Thus, we have
\begin{align}\label{kk}
\langle I'_{\lambda,\mu}(z_{n}),-m^{\star}\rangle+(\xi_{n}(m^{\star})-1)\langle I'_{\lambda,\mu}(z_{n}),z_{n}-m^{\star}\rangle \geq -\frac{1}{n}\|m_{\rho}-z_{n}\|+o(\|m_{\rho}-z_{n}\|).
\end{align}
From $\xi_{n}(m^{\star})(z_{n}-m^{\star})\in \mathbf{\mathbb}{N}_{\lambda,\mu}$ and \eqref{kk}, we obtain
\begin{align*}
-\rho\langle I'_{\lambda,\mu}(z_{n}),\frac{w}{\|w\|}\rangle+(\xi_{n}(m^{\star})-1)\langle I'_{\lambda,\mu}(z_{n})-I'_{\lambda,\mu}(m_{\rho}),z_{n}-m^{\star}\rangle \geq -\frac{1}{n}\|m_{\rho}-z_{n}\|+o(\|m_{\rho}-z_{n}\|).
\end{align*}
So, we get
\begin{align}\label{nnn}
&\langle I'_{\lambda,\mu}(z_{n}),\frac{w}{\|w\|}\rangle\leq\frac{\|m_{\rho}-z_{n}\|}{n\rho}+\frac{o(\|m_{\rho}-z_{n}\|)}{\rho}\\\nonumber
&+\frac{(\xi_{n}(m^{\star})-1)}{\rho}\langle I'_{\lambda,\mu}(z_{n})-I'_{\lambda,\mu}(m_{\rho}),z_{n}-m^{\star}\rangle.
\end{align}
Since$$\|m_{\rho}-z_{n}\|\leq \rho|\xi_{n}(m^{\star})|+|\xi_{n}(m^{\star})-1|\|(u_{n},v_{n})\|$$ and $$\lim_{\rho\rightarrow0}\frac{|\xi_{n}(m^{\star})-1|}{\rho}\leq \|\xi'_{n}(0)\|.$$
If we let $\rho\rightarrow0$ in \eqref{nnn} for fixed $n\in \mathbb{N}$, then by \eqref{abc} we can find a constant $C>0$, independent of $\rho$ such that
$$\langle I'_{\lambda,\mu}(z_{n}),\frac{w}{\|w\|}\rangle\leq\frac{C}{n}\left(1+\|\xi'_{n}(0)\|\right).$$
Thus, we are done once we show that $\|\xi'_{n}(0)\|$ is uniformly bounded. By \eqref{L6}, \eqref{abc} and H\"{o}lder inequality, we have
$$|\langle\xi'_{n}(0),m\rangle|\leq \frac{C_{1}\|m\|}{\left|(p-q)\|(u_{n},v_{n})\|^{p}-2(\alpha+\beta-q)\int_{\Omega}h|u_{n}|^{\alpha}|v_{n}|^{\beta}dx\right|},$$
for some $C_{1}>0$. We only need to show that
$$\left|(p-q)\|(u_{n},v_{n})\|^{p}-2(\alpha+\beta-q)\int_{\Omega}h|u_{n}|^{\alpha}|v_{n}|^{\beta}dx\right|\geq C_{2},$$
for some $C_{2}>0$ and $n$ large enough. We argue by contradiction. Assume that there exists a subsequence $z_{n}$ such that
\begin{align}\label{jj}
(p-q)\|(u_{n},v_{n})\|^{p}-2(\alpha+\beta-q)\int_{\Omega}h|u_{n}|^{\alpha}|v_{n}|^{\beta}dx=o_{n}(1).
\end{align}
By \eqref{jj} and the fact that $z_{n}\in \mathbf{\mathbb}{N}_{\lambda,\mu}$, we have
$$\|(u_{n},v_{n})\|^{p}=\frac{2(\alpha+\beta-q)}{p-q}\int_{\Omega}h|u_{n}|^{\alpha}|v_{n}|^{\beta}dx+o_{n}(1)$$
and $$\|(u_{n},v_{n})\|^{p}=\frac{\alpha+\beta-q}{\alpha+\beta-2}\int_{\Omega}\left(\lambda f|u_{n}|^{q}+\mu g|v_{n}|^{q}\right)dx+o_{n}(1).$$
When n is large enough, by the H\"{o}lder and Sobolev inequalities, we have
\begin{align}\label{M1}
\|(u_{n},v_{n})\|\geq [\frac{p-q}{2(\alpha+\beta-q)}S^{\frac{\alpha+\beta}{p}}]^{\frac{1}{\alpha+\beta-p}}
\end{align}
and
\begin{align*}
\frac{\alpha+\beta-p}{\alpha+\beta-q}\|(u_{n},v_{n})\|^{p}&=\int_{\Omega}\left(\lambda f|u_{n}|^{q}+\mu g|v_{n}|^{q}\right)dx\\
&\leq |\lambda|\|f\|_{L^{q^{\star}}}|\|u_{n}\|^{q}_{L^{\alpha+\beta}}+|\mu|\|g\|_{L^{q^{\star}}}|\|v_{n}\|^{q}_{L^{\alpha+\beta}}\\
&\leq [(|\lambda|\|f\|_{L^{q^{\star}}})^{\frac{p}{p-q}}+(|\mu|\|g\|_{L^{q^{\star}}})^{\frac{p}{p-q}}]^{\frac{p-q}{p}}S^{-\frac{q}{p}}\|(u_{n},v_{n})\|^{q}.
\end{align*}
This implies
\begin{align}\label{M2}
\|(u_{n},v_{n})\|\leq \left(S^{-\frac{q}{p}}\frac{\alpha+\beta-q}{\alpha+\beta-p}\right)^{\frac{1}{p-q}}[(|\lambda|\|f\|_{L^{q^{\star}}})^{\frac{p}{p-q}}+(|\mu|\|g\|_{L^{q^{\star}}})^{\frac{p}{p-q}}]^{\frac{1}{p}}.
\end{align}
By \eqref{M1} and \eqref{M2}, we have
$$[(|\lambda|\|f\|_{L^{q^{\star}}})^{\frac{p}{p-q}}+(|\mu|\|g\|_{L^{q^{\star}}})^{\frac{p}{p-q}}]\geq\left[\frac{p-q}{2(\alpha+\beta-q)}S^{\frac{\alpha+\beta}{p}}\right]^{\frac{p}{\alpha+\beta-p}}\left(S^{-\frac{q}{p}}\frac{\alpha+\beta-q}{\alpha+\beta-p}\right)^{-\frac{p}{p-q}},$$
contradicting the assumption,  that is
$$\langle I'_{\lambda,\mu}(z_{n}),\frac{w}{\|w\|}\rangle\leq\frac{C}{n}.$$This completes the proof of (1).\\
Similarly, by Lemma \ref{M3}, we can prove $(ii)$, we omit the details here.
\end{proof}
\begin{lemma}\label{L14}
If $\{z_{n}\}\subset E$ is a $(PS)_{c}$-sequence for $I_{\lambda,\mu}$, then $\{z_{n}\}$ is bounded in $E$.
\end{lemma}
\begin{proof}
Let $z_{n}=(u_{n},v_{n})\subset E$ be a $(PS)_{c}$-sequence for $I_{\lambda,\mu}$, suppose by contradiction that $\|(u_{n},v_{n})\|\rightarrow+\infty$ as $n\rightarrow+\infty$. Let
$$\widetilde{z_{n}}=(\widetilde{u}_{n},\widetilde{v}_{n}):=\frac{z_{n}}{\|(u_{n},v_{n})\|}=(\frac{u_{n}}{\|(u_{n},v_{n})\|},\frac{v_{n}}{\|(u_{n},v_{n})\|}).$$
We may assume that $\widetilde{z}_{n}\rightharpoonup\widetilde{z}=(\widetilde{u},\widetilde{v})$ in $E$. By the compact embedding theorem, we know $\widetilde{u}_{n}(\cdot,0)\rightarrow\widetilde{u}(\cdot,0)$ and $\widetilde{v}_{n}(\cdot,0)\rightarrow\widetilde{v}(\cdot,0)$ strongly in $L^{r}(\Omega)$ for all $1\leq r<2^{\star}$.
Thus, by H\"{o}lder inequality and Dominated convergence theorem, we have
$$\int_{\Omega}\left(\lambda f|\widetilde{u}_{n}|^{q}+\mu g|\widetilde{v}_{n}|^{q}\right)dx=\int_{\Omega}\left(\lambda f|\widetilde{u}|^{q}+\mu g|\widetilde{v}|^{q}\right)dx+o_{n}(1).$$
Since $\{z_{n}\}$ is a $(PS)_{c}$-sequence for $I_{\lambda,\mu}$ and $\|(u_{n},v_{n})\|\rightarrow+\infty$, we have
\begin{align}\label{L11}
&\frac{1}{p}\int_{Q}\frac{|\widetilde{u}_{n}(x)-\widetilde{u}_{n}(y)|^{p}}{|x-y|^{n+ps}}dxdy+\frac{1}{p}\int_{Q}\frac{|\widetilde{v}_{n}(x)-\widetilde{v}_{n}(y)|^{p}}{|x-y|^{n+ps}}dxdy\\\nonumber
&-\frac{\|(u_{n},v_{n})\|^{q-p}}{q}\int_{\Omega}\left(\lambda f(x)|\widetilde{u}_{n}|^{q}+\mu g(x)|\widetilde{v}_{n}|^{q}\right)dx\\\nonumber
&-\frac{2\|(u_{n},v_{n})\|^{\alpha+\beta-p}}{\alpha+\beta}\int_{\Omega}h(x)|\widetilde{u}_{n}|^{\alpha}|\widetilde{v}_{n}|^{\beta}dx=o_{n}(1)
\end{align}
and
\begin{align}\label{L12}
&\int_{Q}\frac{|\widetilde{u}_{n}(x)-\widetilde{u}_{n}(y)|^{p}}{|x-y|^{n+ps}}dxdy+\int_{Q}\frac{|\widetilde{v}_{n}(x)-\widetilde{v}_{n}(y)|^{p}}{|x-y|^{n+ps}}dxdy\\\nonumber
&-\|(u_{n},v_{n})\|^{q-p}\int_{\Omega}\left(\lambda f(x)|\widetilde{u}_{n}|^{q}+\mu g(x)|\widetilde{v}_{n}|^{q}\right)dx\\\nonumber
&-2\|(u_{n},v_{n})\|^{\alpha+\beta-p}\int_{\Omega}h(x)|\widetilde{u}_{n}|^{\alpha}|\widetilde{v}_{n}|^{\beta}dx=o_{n}(1)
\end{align}
Combining \eqref{L11} and \eqref{L12}, as $n\rightarrow\infty$, we obtain
\begin{align}\label{L13}
&\int_{Q}\frac{|\widetilde{u}_{n}(x)-\widetilde{u}_{n}(y)|^{p}}{|x-y|^{n+ps}}dxdy+\int_{Q}\frac{|\widetilde{v}_{n}(x)-\widetilde{v}_{n}(y)|^{p}}{|x-y|^{n+ps}}dxdy\\\nonumber
&=\frac{p(\alpha+\beta-q)}{q(\alpha+\beta-p)}\|(u_{n},v_{n})\|^{q-p}\int_{\Omega}\left(\lambda f(x)|\widetilde{u}_{n}|^{q}+\mu g(x)|\widetilde{v}_{n}|^{q}\right)dx+o_{n}(1).
\end{align}
Since $1<q<p$ and $\|(u_{n},v_{n})\|\rightarrow+\infty$ as $n\rightarrow+\infty$, \eqref{L13} implies
$$\int_{Q}\frac{|\widetilde{u}_{n}(x)-\widetilde{u}_{n}(y)|^{p}}{|x-y|^{n+ps}}dxdy+\int_{Q}\frac{|\widetilde{v}_{n}(x)-\widetilde{v}_{n}(y)|^{p}}{|x-y|^{n+ps}}dxdy\rightarrow0.$$
Which contradicts the fact that $\|\widetilde{z}_{n}\|=1$ for any $n\geq1$.
\end{proof}
\section{ Local minimization problem }\label{sec5}
Now, we establish the existence of a local minimum for $I_{\lambda,\mu}$ on $\mathbf{\mathbb}{N}^{+}_{\lambda,\mu}$.
\begin{theorem}\label{THE1}
Let $(\lambda,\mu)\in\Theta$, then $I_{\lambda,\mu}$ has a local minimizer $z^{+}$ in $\mathbf{\mathbb}{N}^{+}_{\lambda,\mu}$ satisfying \\
$(i)$ $I_{\lambda,\mu}(z^{+})=\theta_{\lambda,\mu}=\theta^{+}_{\lambda,\mu};$\\
$(ii)$ $z^{+}$ is a nontrivial solution of \eqref{int1}.
\end{theorem}
\begin{proof}
By $(i)$ of Lemma \ref{L10} there exists a minimizing sequence $\{z_{n}\}=\{(u_{n},v_{n})\}$ for $I_{\lambda,\mu}$  in $\mathbf{\mathbb}{N}_{\lambda,\mu}$ such that
\begin{align}\label{L15}
I_{\lambda,\mu}(z_{n})=\theta_{\lambda,\mu}+o_{n}(1) \ and \ I'_{\lambda,\mu}(z_{n})=o_{n}(1) \ in \ E^{-1}.
\end{align}
By Lemma \ref{hhh}, Lemma \ref{L14} and the compact imbedding theorem, we know there is a subsequence, still denoted by $\{z_{n}\}$ and $z^{+}=(u^{+},v^{+})\in E$ such that
\begin{equation*}
\begin{cases}
 u_{n}\rightharpoonup u^{+},v_{n}\rightharpoonup v^{+},$$  &weakly \ \text{in}\ X_{0}^{s}(\Omega),\\
u_{n}\rightarrow u^{+},v_{n}\rightarrow v^{+}, &srongly \ \text{in}\ L^{r}(\Omega)\  for \ all\  1\leq r<2^{\star}.
\end{cases}
\end{equation*}
As $n\rightarrow\infty$, by H\"{o}lder inequality and Dominated convergence theorem, we obtain
\begin{align}\label{L16}
\int_{\Omega}\left(\lambda f|u_{n}|^{q}+\mu g|v_{n}|^{q}\right)dx=\int_{\Omega}\left(\lambda f|u^{+}|^{q}+\mu g|v^{+}|^{q}\right)dx+o_{n}(1)
\end{align}
and
\begin{align}\label{L17}
\int_{\Omega}h|u_{n}|^{\alpha}|v_{n}|^{\beta}dx=\int_{\Omega}h|u^{+}|^{\alpha}|v^{+}|^{\beta}dx+o_{n}(1).
\end{align}
First, we claim that $\int_{\Omega}\left(\lambda f|u^{+}|^{q}+\mu g|v^{+}|^{q}\right)dx\neq0$, we argue by contradiction, then we have $\int_{\Omega}\left(\lambda f|u_{n}|^{q}+\mu g|v_{n}|^{q}\right)dx\rightarrow0$ as $n\rightarrow\infty$. Thus$$\|(u_{n},v_{n})\|^{p}=2\int_{\Omega}h|u_{n}|^{\alpha}|v_{n}|^{\beta}dx+o_{n}(1)$$
and
\begin{align*}
&I_{\lambda,\mu}(z_{n})=\frac{1}{p}\|(u_{n},v_{n})\|^{p}-
\frac{2}{\alpha+\beta}\int_{\Omega}h(x)|u_{n}|^{\alpha}|v_{n}|^{\beta}dx+o_{n}(1)\\
&=\left(\frac{1}{p}-\frac{1}{\alpha+\beta}\right)\|(u_{n},v_{n})\|^{2}+o_{n}(1).
\end{align*}
This contradicts $I_{\lambda,\mu}(z_{n})\rightarrow\theta_{\lambda,\mu}<0$ as $n\rightarrow\infty$.\\
Now, we claim $z^{+}$ is a nontrivial solution of \eqref{int1}. From \eqref{L15}, \eqref{L16} and \eqref{L17}, we know $z^{+}$ is a weak solution of \eqref{int1}. From $z_{n}\in \mathbf{\mathbb}{N}_{\lambda,\mu}$, we have
\begin{align}\label{L20}
I_{\lambda,\mu}(z_{n})=\frac{\alpha+\beta-p}{p(\alpha+\beta)}\|(u_{n},v_{n})\|^{p}-\frac{\alpha+\beta-q}{q(\alpha+\beta)}\int_{\Omega}\left(\lambda f|u_{n}|^{q}+\mu g|v_{n}|^{q}\right)dx.
\end{align}
That is
\begin{align}\label{L18}
\int_{\Omega}\left(\lambda f|u_{n}|^{q}+\mu g|v_{n}|^{q}\right)dx=\frac{q(\alpha+\beta-p)}{p(\alpha+\beta-q)}\|(u_{n},v_{n})\|^{p}-\frac{q(\alpha+\beta)}{\alpha+\beta-q}I_{\lambda,\mu}(z_{n}).
\end{align}
Let $n\rightarrow\infty$ in \eqref{L18}, by \eqref{L15}, \eqref{L16} and $\theta_{\lambda,\mu}<0$, we have$$\int_{\Omega}\left(\lambda f|u^{+}|^{q}+\mu g|v^{+}|^{q}\right)dx\geq -\frac{q(\alpha+\beta)}{\alpha+\beta-q}\theta_{\lambda,\mu}>0.$$
Therefore, $z^{+}\in \mathbf{\mathbb}{N}_{\lambda,\mu}$ is a nontrival solution of \eqref{int1}.
Next, we show that $z_{n}\rightarrow z^{+}$ strongly in $E$ and $I_{\lambda,\mu}(z^{+})=\theta_{\lambda,\mu}$. Since $z^{+}\in \mathbf{\mathbb}{N}_{\lambda,\mu}$, then by \eqref{L20}, we obtain
\begin{align}
\theta_{\lambda,\mu}\leq I_{\lambda,\mu}(z^{+})&=\frac{\alpha+\beta-p}{p(\alpha+\beta)}\|(u^{+},v^{+})\|^{p}-\frac{\alpha+\beta-q}{q(\alpha+\beta)}\int_{\Omega}\left(\lambda f|u^{+}|^{q}+\mu g|v^{+}|^{q}\right)dx\\\nonumber
&\leq \liminf\limits_{n\rightarrow \infty}\left(\frac{\alpha+\beta-p}{p(\alpha+\beta)}\|(u_{n},v_{n})\|^{p}-\frac{\alpha+\beta-q}{q(\alpha+\beta)}\int_{\Omega}\left(\lambda f|u_{n}|^{q}+\mu g|v_{n}|^{q}\right)dx\right)\\\nonumber
&\leq \lim\limits_{n\rightarrow \infty}\left(\frac{\alpha+\beta-p}{p(\alpha+\beta)}\|(u_{n},v_{n})\|^{p}-\frac{\alpha+\beta-q}{q(\alpha+\beta)}\int_{\Omega}\left(\lambda f|u_{n}|^{q}+\mu g|v_{n}|^{q}\right)dx\right)\\\nonumber
&\leq \lim\limits_{n\rightarrow \infty}I_{\lambda,\mu}(z_{n})=\theta_{\lambda,\mu}.
\end{align}
This inplies that $I_{\lambda,\mu}(z^{+})=\theta_{\lambda,\mu}$ and $\lim\limits_{n\rightarrow\infty}\|(u_{n},v_{n})\|^{p}=\|(u^{+},v^{+})\|^{p}$. Hence $z_{n}\rightarrow z^{+}$ srongly in $E$.\\
Finially, we claim that $z^{+}\in \mathbf{\mathbb}{N}^{+}_{\lambda,\mu}$. Assume by contradiction that $z^{+}\in \mathbf{\mathbb}{N}^{-}_{\lambda,\mu}$, then by Lemma \ref{L21}, there exist unique $t_{1}^{+}$ and $t_{1}^{-}$, such that $t_{1}^{+}(z^{+})\in \mathbf{\mathbb}{N}^{+}_{\lambda,\mu},\ t_{1}^{-}(z^{+})\in \mathbf{\mathbb}{N}^{-}_{\lambda,\mu}$. In particular, we have $t_{1}^{+}<t_{1}^{-}=1$.
Since $$\frac{d}{dt}I_{\lambda,\mu}(t_{1}^{+}z^{+})=0 \ and \  \frac{d^{2}}{dt^{2}}I_{\lambda,\mu}(t_{1}^{+}z^{+})>0,$$ there exists $t_{1}^{+}<t^{\star}<t_{1}^{-}$ such that $I_{\lambda,\mu}(t_{1}^{+}z^{+})<I_{\lambda,\mu}(t^{\star}z^{+})$. By Lemma \ref{L21}, we have $$I_{\lambda,\mu}(t_{1}^{+}z^{+})<I_{\lambda,\mu}(t^{\star}z^{+})\leq I_{\lambda,\mu}(t_{1}^{-}z^{+})=I_{\lambda,\mu}(z^{+}),$$ a contraction. Since $I_{\lambda,\mu}(z^{+})=I_{\lambda,\mu}(|u^{+}|,|v^{+}|)$ and $(|u^{+}|,|v^{+}|)\in \mathbf{\mathbb}{N}_{\lambda,\mu}$, by Lemma \ref{L22}, we obtain that $z^{+}$ is a nontrivial solution of \eqref{int1}.
\end{proof}

Next, we establish the existence of a local minimum for $I_{\lambda,\mu}$ on $\mathbf{\mathbb}{N}^{-}_{\lambda,\mu}$.
\begin{theorem}\label{THE2}
Let $(\lambda,\mu)\in\Psi$, then $I_{\lambda,\mu}$ has a local minimizer $z^{-}$ in $\mathbf{\mathbb}{N}^{-}_{\lambda,\mu}$ satisfying \\
$(i)$ $I_{\lambda,\mu}(z^{-})=\theta^{-}_{\lambda,\mu};$\\
$(ii)$ $z^{-}$ is a nontrivial solution of \eqref{int1}.
\end{theorem}
\begin{proof}
By $(ii)$ of Lemma \ref{L10} there exists a minimizing sequence $\{z_{n}\}=\{(u_{n},v_{n})\}$ for $I_{\lambda,\mu}$  in $\mathbf{\mathbb}{N}^{-}_{\lambda,\mu}$ such that
\begin{align*}
I_{\lambda,\mu}(z_{n})=\theta^{-}_{\lambda,\mu}+o_{n}(1) \ and \ I'_{\lambda,\mu}(z_{n})=o_{n}(1) \ in \ E^{-1}.
\end{align*}
By Lemma \ref{hhh} $(iii)$, Lemma \ref{L14} and the compact imbedding theorem, we know there is a subsequence, still denoted by $\{z_{n}\}$ and $z^{-}=(u^{-},v^{-})\in \mathbf{\mathbb}{N}^{-}_{\lambda,\mu}$ such that
\begin{equation*}
\begin{cases}
 u_{n}\rightharpoonup u^{-},v_{n}\rightharpoonup v^{-},$$  &weakly \ \text{in} X_{0}^{s}(\Omega),\\
u_{n}\rightarrow u^{-},v_{n}\rightarrow v^{-}, &srongly \ \text{in} L^{r}(\Omega)\  for \ all\  1\leq r<2^{\star}.
\end{cases}
\end{equation*}
As $n\rightarrow\infty$, by H\"{o}lder inequality and Dominated convergence theorem, we obtain
\begin{align*}
\int_{\Omega}\left(\lambda f|u_{n}|^{q}+\mu g|v_{n}|^{q}\right)dx=\int_{\Omega}\left(\lambda f|u^{-}|^{q}+\mu g|v^{-}|^{q}\right)dx+o_{n}(1)
\end{align*}
and
\begin{align*}
\int_{\Omega}h|u_{n}|^{\alpha}|v_{n}|^{\beta}dx=\int_{\Omega}h|u^{-}|^{\alpha}|v^{-}|^{\beta}dx+o_{n}(1).
\end{align*}
First, we claim that $\int_{\Omega}\left(\lambda f|u^{-}|^{q}+\mu g|v^{-}|^{q}\right)dx\neq0$, suppose by contradiction, then we have $\int_{\Omega}\left(\lambda f|u_{n}|^{q}+\mu g|v_{n}|^{q}\right)dx\rightarrow0$ as $n\rightarrow\infty$. Thus$$\|(u_{n},v_{n})\|^{p}=2\int_{\Omega}h|u_{n}|^{\alpha}|w_{2,n}|^{\beta}dx+o_{n}(1)$$
and
\begin{align*}
&I_{\lambda,\mu}(z_{n})=\frac{1}{p}\|(u_{n},v_{n})\|^{2}-
\frac{2}{\alpha+\beta}\int_{\Omega}h(x)|u_{n}|^{\alpha}|v_{n}|^{\beta}dx+o_{n}(1)\\
&=\left(\frac{1}{p}-\frac{1}{\alpha+\beta}\right)\|(u_{n},v_{n})\|^{2}+o_{n}(1).
\end{align*}
This contradicts $I_{\lambda,\mu}(z_{n})\rightarrow\theta_{\lambda,\mu}<0$ as $n\rightarrow\infty$.\\
Now, we prove that $z_{n}\rightarrow z^{-}$ strongly in $E$. Othercase, we have
\begin{align*}
&\|(u^{-},v^{-})\|^{p}-\int_{\Omega}\left(\lambda f|u^{-}|^{q}+\mu g|v^{-}|^{q}\right)dx-2\int_{\Omega}h|u^{-}|^{\alpha}|v^{-}|^{\beta}dx\\
&\leq\liminf\limits_{n\rightarrow \infty}\left(\|(u_{n},v_{n})\|^{p}-\int_{\Omega}\left(\lambda f|u_{n}|^{q}+\mu g|v_{n}|^{q}\right)dx-2\int_{\Omega}h|u_{n}|^{\alpha}|v_{n}|^{\beta}dx\right)\\
&\leq\lim\limits_{n\rightarrow \infty}\left(\|(u_{n},v_{n})\|^{p}-\int_{\Omega}\left(\lambda f|u_{n}|^{q}+\mu g|v_{n}|^{q}\right)dx-2\int_{\Omega}h|u_{n}|^{\alpha}|v_{n}|^{\beta}dx\right)=0.
\end{align*}
Which contradicts $z^{-}\in \mathbf{\mathbb}{N}^{-}_{\lambda,\mu}$. Hence $z_{n}\rightarrow z^{-}$ strongly in $E$.
This implies $$I_{\lambda,\mu}(z_{n})\rightarrow I_{\lambda,\mu}(z^{-})=\theta^{-}_{\lambda,\mu} \ as\ n\rightarrow +\infty.$$
Since $I_{\lambda,\mu}(z^{-})= I_{\lambda,\mu}(|u^{-},v^{-}|)$ and $|u^{-},v^{-}|\in \mathbf{\mathbb}{N}^{-}_{\lambda,\mu}$, by Lemma \ref{L22}, we have $z^{-}$ is a nontrivial solution of \eqref{int1}.
\end{proof}
\section{Proof of Theorem \ref{TH1} and Theorem \ref{TH2}}\label{sec6}
Now, we complete the proof of Theorem \ref{TH1} and Theorem \ref{TH2}.
\begin{proof}
 For $(\lambda,\mu)\in\Theta$, by Theorem \ref{THE1}, system \eqref{int1} admits at least one nontrivial solution $z^{+}\in\mathbf{\mathbb}{N}^{+}_{\lambda,\mu}$. By Theorem \ref{THE1} and Theorem \ref{THE2}, we obtain that for $(\lambda,\mu)\in\Psi$, system \eqref{int1} admits at least two nontrivial solutions  $z^{+}$ and $z^{-}$ such that $z^{+}\in\mathbf{\mathbb}{N}^{+}_{\lambda,\mu}$,  $z^{-}\in\mathbf{\mathbb}{N}^{-}_{\lambda,\mu}$. Since $\mathbf{\mathbb}{N}^{+}_{\lambda,\mu}\cap\mathbf{\mathbb}{N}^{+}_{\lambda,\mu}=\emptyset$, then $z^{+}$ and $z^{-}$ are distinct solutions of syetem \eqref{int1}.
\end{proof}

\end{document}